\magnification =1200

\documentstyle{amsppt}
\NoBlackBoxes

\pageheight{9truein}
\pagewidth{6.5truein}

\newcount\sectct \sectct =0
\newcount\thmct \thmct =0

\def\section#1{\global\thmct =0 \global\advance \sectct by 1 \vskip 4ex
		\centerline{\uppercase\expandafter{\romannumeral\sectct.}
		\uppercase{ #1}} \vskip 2ex}

\def\prclaim#1{\global \advance \thmct by 1 \vskip 2ex
        \noindent{\bf #1\ \the\sectct.\the\thmct.}
                \unskip \bgroup \it \unskip}
\def\endprclaim{{\egroup} \vskip 0ex}
\def\dmo#1{{\global \advance \thmct by 1 \vskip 1ex
        \noindent{\bf #1\ \the\sectct.\the\thmct.}}}
\def\label#1{{\unskip\global\edef#1{\the\sectct.\the\thmct\unskip}}}

\define\perm{P_2(M)}
\def\Sfg{S(f, g)}
\def\ijkm{1\leq i<j<k\leq m}

\def\xipjqkr{x_{ip}x_{jq}x_{kr}}

\def\pqrn{1\leq p<q<r\leq n}
\def\dpqrn{\text{$p, q, r$ distinct}}
\def\dijkm{\text{$i, j, k$ distinct}}

\document

\topmatter
\title
Permanental Ideals
\endtitle

\author
Reinhard C. Laubenbacher and Irena Swanson
\endauthor
\address
Reinhard C. Laubenbacher and Irena Swanson, Department of Mathematics,
New Mexico State University,
Las Cruces, NM 88003
\endaddress
\email
reinhard\@nmsu.edu iswanson\@nmsu.edu
\endemail


\thanks
The authors thank Bernd Sturmfels for drawing their attention to
the subject of this paper, the Mathematical Sciences Research Institute, Berkeley,
for its hospitality
while part of this paper was being written,
Eberhard Becker
for drawing their attention to Michael Niermann's thesis, and Luchezar Avramov for
helpful comments.
The second author also thanks the NSF for partial support.
\endthanks

\abstract
The principal result is a primary decomposition of ideals
generated by the $(2\times 2)$-subpermanents of a generic matrix.
These permanental ideals almost always have embedded components
and their minimal primes are of three distinct heights.
Thus the permanental ideals are almost never Cohen-Macaulay,
in contrast with the determinantal ideals.
\endabstract

\endtopmatter
\document

\section{Introduction.}

Permanents were introduced by Cauchy and
Binet at the beginning of the nineteenth century as a special type
of alternating symmetric function.  Later they were studied by
I. Schur as a special type of what are now known as Schur functions.
Permanents have since found applications in combinatorics,
probability theory, and, more recently, invariant theory \cite{HK}.
A good survey of the theory of permanents is \cite{Mi}.

In contrast to determinantal ideals,
permanental ideals have not received much attention to date.
M.\ Niermann, in his Ph.D.\ thesis \cite{N},
calculates the radicals and real radicals
of $2\times 2$ generic permanental ideals.
Niermann's motivation,
as well as ours,
came from work by D. Eisenbud and B. Sturmfels on binomial
ideals \cite{ES}.

The {\it permanent} of an $(n\times n)$ square matrix $M=(a_{ij})$ 
is defined as
$$
\text{perm}(M)=\sum_{\sigma \in S_n}a_{1\sigma (1)}a_{2\sigma (2)}\cdots
                  a_{n\sigma (n)}.
$$
Thus, the permanent differs from the determinant only in the lack of
minus signs in the expansion.

In this paper we work with generic matrices $M$.
More precisely,
let $F$ be a field,
$m, n, r$ positive integers,
and $x_{ij}$ variables over $F$
with $1 \leq i\leq m$ and $1 \leq j \leq n$.
Let $R=F[x_{ij}|1\leq i\leq m,1\leq j\leq n]$ be the polynomial ring,
and let $M$ be the $(m\times n)$-matrix whose $(i,j)$-entry is $x_{ij}$.
Then let $P_r(M)$ be the ideal of $R$
generated by all the $(r\times r)$-subpermanents of $M$.
If the field $F$ has characteristic $2$,
the permanental ideals equal determinantal ideals,
which are relatively well-understood \cite{BH, BV}.
Thus for the rest of this paper
we assume that the characteristic of $F$ is different from $2$.

In this paper we study the binomial ideals $\perm$.
Their properties turn out to be
very different from those
of the much better understood determinantal ideals.
For example,
we prove that the generating permanents of $\perm$
are not a Gr\"obner basis of $\perm$
in any diagonal order,
whereas Caniglia, Guccione and Guccione \cite{CGG}
and, independently, Sturmfels \cite{S}
proved that the generating $(2\times 2)$-minors of $M$
are a Gr\"obner basis of the ideal they generate in any diagonal order.
Also, in contrast to determinantal ideals, 
we prove that if $m$ and $n$ are both at least 3,
then $\perm$ is not a radical ideal,
is not Cohen-Macaulay,
and there are minimal primes of distinct heights over it.
Furthermore,
we prove that
the radical of $\perm$ requires generators of degree higher than
those of $\perm$ itself.
If $m, n \ge 4$,
$\perm$ is not integrally closed.
We explicitly calculate a primary decomposition,
the unmixed parts of $\perm$ of all possible dimensions,
and the radical.
Radicals and real radicals
in characteristic $0$
were also calculated independently by Niermann \cite{N}.

Thus, the permanental ideals are another case of ideals
for which the Gr\"obner bases,
irreducibility,
primary decompositions,
Cohen-Macaulayness,
and integral closedness
depend on the characteristic of the underlying field:
we have ``good'' properties in characteristic $2$
versus very different properties in all other characteristics.
Moreover,
for all the properties that we study,
the results are independent of the characteristic as
long as the characteristic is not $2$!

\section{Monomials in $\perm$.}

In this section we prove that, unlike the determinantal ideals,
the ideals $\perm$ contain many monomials.
The proof depends heavily on the characteristic
of the field not being $2$.
There are two basic types of monomials in $\perm$,
treated in the two lemmas below.

\prclaim{Lemma}
\label{\twoone}
The ideal $\perm$ contains all products of three entries of $M$,
taken from three distinct columns and two distinct rows,
or from two distinct columns and three distinct rows.
\endprclaim

\demo{Proof}
We only prove the first part.  Without changing $P_2(M)$ we may permute rows
and columns of $M$, and so we assume that the 
three entries are contained in a submatrix
$$
\pmatrix a&b&c\\x&y&z\endpmatrix 
$$
of $M$.
Then $\perm$ contains the elements $ay+bx$ and $bz+cy$.  Therefore,
$$
c(ay+bx)-a(bz+cy) = cbx-abz = b(cx-az)\in \perm.
$$
Since $cx+az\in \perm$ and 2 is a unit in $R$, 
we obtain that $bcx$ and $abz$ are in $\perm$.
A similar argument shows that all other products
as in the statement are in $\perm$ as well.
This proves the lemma.
\enddemo

\prclaim{Lemma}
\label{\squarethree}
If $m, n \ge 3$,
$\perm$ contains all monomials of the form
$x_{i_1j_1}^{e_1}x_{i_2j_2}^{e_2}x_{i_3j_3}^{e_3}$
with distinct $i_1, i_2, i_3$,
distinct $j_1, j_2, j_3$,
and where $e_1$, $e_2$, and $e_3$ are positive integers
which sum to $4$.
\endprclaim

\demo{Proof}
Permuting rows and columns as before, we
consider the submatrix
$$
\pmatrix a&b&c\\x&y&z\\ u&v&w \endpmatrix .
$$
It suffices to prove that $ayw^2$ lies in $\perm$.
But
$$
ayw^2 = aw(yw+vz) - awvz.
$$
As $yw+vz \in \perm$ by definition
and $awvz \in \perm$ by Lemma~\twoone,
we are done.
\enddemo

\section{A Gr\"obner Basis.}

We compute a Gr\"obner basis for $\perm$ with respect
to any lexicographic diagonal ordering of monomials.
Recall that a monomial order on the $x_{ij}$ is {\it diagonal}
if for any square submatrix of $M$,
the leading term of the permanent (or of the determinant)
of that submatrix is the product of the entries on the main diagonal.
An example of such an order is the lexicographic order defined by:
$$
x_{ij}<x_{kl} \text{ if and only if }l>j\text{ or }l = j\text{ and }k>i.
$$
Throughout this section we use an arbitrary lexicographic diagonal
ordering.

\prclaim{Theorem}
\label{\gb}
The following collection $G$ of polynomials
is a minimal reduced Gr\"obner basis for $\perm$,
with respect to any diagonal ordering:
\roster
\item
The subpermanents $x_{ij}x_{kl}+x_{kj}x_{il}, i<k, j<l$;
\item
$x_{i_1j_1}x_{i_1j_2}x_{i_2j_3}, i_1>i_2, j_1<j_2<j_3$;
\item
$x_{i_1j_1}x_{i_2j_2}x_{i_2j_3}, i_1>i_2, j_1<j_2<j_3$;
\item
$x_{i_1j_1}x_{i_2j_1}x_{i_3j_2}, i_1<i_2<i_3, j_1>j_2$;
\item
$x_{i_1j_1}x_{i_2j_2}x_{i_3j_2}, i_1<i_2<i_3, j_1>j_2$;
\item
$x_{i_1j_1}^ax_{i_2j_2}^bx_{i_3j_3}^c, i_1<i_2<i_3, j_1>j_2>j_3,
abc =2$.
\endroster
\endprclaim

These monomials can also be described as follows.
Monomials of type (6)
are in one-to-three correspondence with the $3$ by $3$ submatrices of $M$:
for each $3$ by $3$ submatrix of $M$,
a monomial of type (6) is the product of the entries
on its anti-diagonal,
with one of the entries taken to the second power.
Pictorially,
in a given $3$ by $3$ submatrix of $M$,
a monomial of type (6)
is the product of the following entries marked $\circ$:
$$
\pmatrix \ & \ & \circ^\cdot \cr
		\  & \circ^\cdot & \  \cr
		\circ^\cdot & \ & \  \cr
\endpmatrix
$$
(The superscripts $\ ^\cdot$ are a reminder that
one of the entries is raised to the second power.)
Similarly,
the monomials of types (2) through (5) are products of
the following entries, marked $\circ$,
of appropriately sized submatrices:
$$
\matrix
\pmatrix \ & \ & \circ \cr
		\circ & \circ & \  \cr
\endpmatrix
&
\hskip 3em
\pmatrix \ & \circ & \circ \cr
		\circ & \ & \  \cr
\endpmatrix
\hskip 3em
&
\pmatrix \ & \circ \cr
		\ & \circ \cr
		\circ & \  \cr
\endpmatrix
\hskip 3em
&
\pmatrix \ & \circ \cr
		\circ & \  \cr
		\circ & \  \cr
\endpmatrix
\cr
\cr
\text{type (2)} & \text{type (3)} & \hskip -3em 
\text{type (4)} & \text{type (5)} \cr
\endmatrix
$$

With this pictorial representation of the monomials in $\perm$
it is easy to count the number of elements
in the minimal Gr\"obner basis in the theorem,
which is
$$
{m \choose 2} {n \choose 2}
+ 2 {m \choose 2}{n \choose 3}
+ 2  {n \choose 2} {m\choose 3}
+ 3 {m \choose 3} {n \choose 3}.
$$

\demo{Proof of Theorem \gb}
By Lemmas~\twoone\ and \squarethree,
$G$ is contained in $\perm$.
As all elements of the form (1) generate $\perm$,
$G$ certainly is a generating set.
It is easy to see that the set is reduced and minimal.
Thus it is sufficient to show that all the $S$-polynomials $\Sfg$ of 
pairs of elements of $G$ reduce to zero with respect to $G$.
As the $S$-polynomials of two monomials always reduce to zero,
it suffices to prove that $\Sfg$ reduces to zero
whenever $f$ is of type (1).

First consider 
$$
\align
f&= x_{ij}x_{kl}+x_{kj}x_{il},\\
g&= x_{pq}x_{rs}+x_{rq}x_{ps},
\endalign
$$
with $i<k, j<l$ and $p<r, q<s$.  For both $f$ and $g$ the first summand
is the leading term.  If $\text{in}(f) = \text{in}(g)$, then $f = g$,
so $\Sfg =0$.
If $\text{in}(f)$ and $\text{in}(g)$ have no factors in common,
then by standard Gr\"obner basis theory $\Sfg$ reduces to zero.

So we may assume that $\text{in}(f)$ and $\text{in}(g)$ have
exactly one variable in common.  Suppose first that $x_{ij}= x_{pq}$.
Then
$$
\align
\Sfg&= x_{rs}f-x_{kl}g\\
      &= x_{rs}x_{kj}x_{il}-x_{kl}x_{rq}x_{ps}\\
	  &= x_{rs}x_{kj}x_{il}-x_{kl}x_{rj}x_{is}.
\endalign
$$
If $r = k$, then 
$$
\Sfg = x_{rs}x_{rj}x_{il}-x_{rl}x_{rj}x_{is}.
$$
The first summand is an element of $G$
of type (3) and the second is of type (2),
thus $\Sfg$ reduces to $0$.
A similar argument also shows that $\Sfg$ reduces to $0$ in case
$j = s$,
but there the reduction is with respect to elements of $G$
of types (4) and (5).
So we may now assume that $r \not = k$, $j \not = s$,
and, by switching,
we may furthermore assume without loss of generality that $k < r$.
Then the leading term of $\Sfg$ is $x_{rs}x_{kj}x_{il}$, and
after subtracting 
$(x_{kj}x_{rs}+x_{rj}x_{ks})x_{il}$,
$\Sfg$ reduces to
$-x_{rj}(x_{ks}x_{il}+x_{kl}x_{is})$, which in turn reduces to $0$
with respect to
the element $x_{ks}x_{il}+x_{kl}x_{is}$ of $G$.

Now suppose that the common factor of $\text{in}(f)$ and $\text{in}(g)$
is $x_{ij} = x_{rs}$.  Then $i = r, j = s$, and
$$
\Sfg = x_{pq}f-x_{kl}g = x_{pq}x_{kj}x_{il}-x_{kl}x_{iq}x_{pj}.
$$
Since $p < r = i$ and $q < s = j < l$,
we may subtract from $\Sfg$
the polynomial $x_{kj} (x_{pq}x_{il} + x_{pl}x_{iq})$
(a monomial multiple of an element of type (1))
to get
$$
-x_{kj}x_{pl}x_{iq} -x_{kl}x_{iq}x_{pj}
= - x_{iq}( x_{kj}x_{pl} + x_{kl}x_{pj}),
$$
whence $\Sfg$ reduces to zero also in this case.
The case $x_{kl} = x_{pq}$ is done similarly.

Finally
(for the case $f$, $g$ both permanents)
we discuss the case when the common factor of the initial terms
is $x_{kl}= x_{rs}$.
Then $k = r, l = s$, and
$$
\Sfg = x_{pq}x_{kj}x_{il}-x_{ij}x_{kq}x_{pl}.
$$
First suppose that $i = p$.
Then without loss of generality $j < q$
and
$\Sfg = x_{iq}x_{kj}x_{il}-x_{ij}x_{kq}x_{il}$.
The first summand is an element of $G$ of type (3),
and the second reduces to zero with respect to $G$, as we obtain an element of $G$
of type (3) by first adding
$(x_{ij}x_{kq} + x_{iq}x_{kj}) x_{il}$ to it
to get $x_{iq}x_{kj} x_{il}$.
Similarly,
if $j = q$,
$\Sfg$ reduces to $0$ with respect to elements of $G$
of type (4).

So it remains to consider the cases $i < p$, $j \not = q$
(and still $k = r$, $l = s$).
Then $i<p<r = k$ and $j<l$.
There are two cases to consider here: $j < q$ and $j > q$.
In both cases $\Sfg$ reduces to $0$
with respect to elements of $G$ of type (1).
This shows that the $S$-polynomial of any two subpermanents reduces
to zero.

Next we prove that the $S$-polynomial of any subpermanent
and any element of $G$ of types (2) - (6) reduces to $0$.
Let $f$ be as before, and let $g$ be a monomial in $G$.
Note that for all $g$ of types (2) through (6),
$\text{in}(f)$ does not divide $g$.
If $\text{in}(f)$ and $g$ have no common factors,
then obviously $\Sfg$ reduces to $0$.
So it suffices to consider only the cases when
$\text{in}(f)$ and $g$ have exactly one variable in common.

First assume that $g$ is of type (2) or (3).
By similarity this will also cover the cases when $g$ is
of types (4) or (5).
So let
$g=x_{i_1j_1}x_{i_1j_2}x_{i_2j_3}$ be of type (2), with $i_1>i_2$ and
$j_1<j_2<j_3$, and $f=x_{ij}x_{kl}+x_{kj}x_{il}$, with $i<k,j<l$.  
Suppose first that $i_1=i,j_1=j$, and $x_{kl}$ does not divide $g$.
Then 
$$
S(f,g)=x_{ij_2}x_{i_2j_3}x_{kj}x_{il}.
$$
If $l=j_2$, then $S(f,g)$ is a monomial of type (6)
so it reduces to $0$.
If $l>j_3$, then $S(f,g)$ reduces to $x_{kj}x_{ij_2}x_{ij_3}x_{i_2l}$
since it is divisible by the leading term of a subpermanent.
This monomial, in turn, reduces to zero, since it is divisible by a 
monomial of type (3).
Likewise, if $l\le j_3$ and $l$ is different from $j_2$,
then $S(f,g)$ is a multiple of a monomial of type (3).
Thus, in all cases $S(f,g)$ reduces to $0$.

Next suppose that $i=i_1,j=j_2$.  Again, if $l>j_3$, then 
$S(f,g)$ is divisible by the leading term of a permanent, and the
resulting reduced monomial is divisible by a monomial of type (2).
If $l=j_3$, then $S(f,g)$ is a multiple of a monomial of type (4),
and if $l<j_3$, then $S(f,g)$ is a multiple of a monomial of either
type (2) or type (6).
Finally, suppose that $i=i_2,j=j_3$.  Then $S(f,g)$ is a 
multiple of a monomial of type (3).
Similar computations show that $S(f,g)$ reduces to zero, if the
common factor of $f$ and $g$ is $x_{kl}$, and for $g$ of type (3).

Now consider a monomial $g$ of type (6), say 
$g=x_{i_1j_1}^2x_{i_2j_2}x_{i_3j_3}$, with $i_1<i_2<i_3$ and 
$j_1>j_2>j_3$.  If $i=i_1,j=j_1$, then $S(f,g)$ is a multiple of
a monomial of type (3).  Now suppose that $i=i_2, j=j_2$.  
If $l>j_1$, respectively $l=j_1$, then $S(f,g$ is divisible by
a monomial of type (6), respectively of type (5).  If $l>j_1$, then
$S(f,g)$ is divisible by the leading term of a permanent, and the
reduced monomial is divisible by a monomial of type (2).
Similar computations show that in all other cases $S(f,g)$ reduces to zero
with respect to $G$.
This completes the proof of Theorem 3.1.
\enddemo

This theorem is in contrast to the case of determinantal ideals
for which the determinants themselves form the reduced Gr\"obner basis
in any diagonal order.
The determinantal result was proved
independently by Caniglia, Guccione and Guccione \cite{CGG}
and by Sturmfels \cite{S}.

\demo{Remark}
It is an easy consequence of this theorem that, if $m, n \ge 4$,
then $\perm$ is not integrally closed.
Lemma \squarethree\ %
shows that $x_{41}x_{32}x_{23}^2$,
$x_{41}x_{32}x_{14}^2$ are both contained in $\perm$.
Then the element
$x_{41}x_{32}x_{23}x_{14}$ is in the integral closure of $\perm$
(as its square equals the product of the first two monomials).
However,
the previous theorem says that
$x_{41}x_{32}x_{23}x_{14}$ is not a multiple of any initial
term of any element of the minimal Gr\"obner basis of $\perm$,
thus it is not in $\perm$.
\enddemo

A slight modification of the proof of Theorem~\gb\ also gives:

\prclaim{Theorem}
\label{\gbrad}
Assume that $m, n \ge 3$.
Let $I$ be the ideal of $R$ generated by elements of $\perm$
and all the products of three entries of $M$
taken from three distinct rows and three distinct columns.
The following collection of polynomials
is a minimal reduced Gr\"obner basis for $I$
with respect to any diagonal ordering:
\roster
\item
The subpermanents $x_{ij}x_{kl}+x_{kj}x_{il}, i<k, j<l$;
\item
$x_{i_1j_1}x_{i_1j_2}x_{i_2j_3}, i_1>i_2, j_1<j_2<j_3$;
\item
$x_{i_1j_1}x_{i_2j_2}x_{i_2j_3}, i_1>i_2, j_1<j_2<j_3$;
\item
$x_{i_1j_1}x_{i_2j_1}x_{i_3j_2}, i_1<i_2<i_3, j_1>j_2$;
\item
$x_{i_1j_1}x_{i_2j_2}x_{i_3j_2}, i_1<i_2<i_3, j_1>j_2$;
\item
$x_{i_1j_1}x_{i_2j_2}x_{i_3j_3}, i_1<i_2<i_3, j_1>j_2>j_3$.
\endroster
\endprclaim

We will show in Theorem 5.4
that the ideal $I$ is equal to the radical of $\perm$
(when $m, n \ge 3$).
With reasoning as in the remark above
we get the following information about the module ${I /\perm}$:

\prclaim{Theorem}
\label{\thmfinite}
Let $m, n \ge 3$.
Let $I$ be the ideal 
of Theorem~\gbrad.
Then
\roster
\item
The elements of type (6) as in Theorem~\gbrad\ %
are a set of minimal generators of the module ${I \over \perm}$.
Thus the number of minimal generators of
${I \over \perm}$ is
${m \choose 3} {n \choose 3}$.
\item
The module ${I \over \perm}$ has finite length as an $F$-vector space.
In fact,
${I \over \perm}$
is minimally generated over $F$ by the
products of the anti-diagonal entries of each $(i\times i)$ 
submatrix of $M$
as $i$ varies from $3$ to $\min\{m,n\}$.
Thus
$$
\eqalignno{
\hbox{length\,}
\left({I \over \perm}\right)
&=
\sum_{i \ge 3} \hbox{(number of $(i\times i)$ submatrices of $M$)} \cr
&=
\sum_{i \ge 3} {m \choose i} {n \choose i} .\cr
}
$$
\endroster
\endprclaim

\section{Minimal Primes.}

We compute the primes of $R$ minimal over $\perm$.
Our computation does not rely on any monomial ordering of the variables,
thus in the proofs we may and do
take transposes of $M$
and we do permute the columns and rows of $M$ as needed.

\prclaim{Theorem} 
\label{\minp}
Let $m, n\geq 2$.
Each of the prime ideals $P$ of $R$ minimal over
$\perm$ is one of the following:
\roster
\item
If $n\geq 3$, then $P$ is generated by all the indeterminates in 
$m-1$ of the rows of $M$;
\item
if $m\geq 3$,
then $P$ is generated by all the indeterminates in $n-1$ of
the columns of $M$;
\item
$P$ is generated by the permanent of one $(2\times 2)$-submatrix 
of $M$ and all the entries of $M$ outside of this submatrix.
\endroster
Moreover,
each of the primes in (1), (2) and (3) is minimal over $\perm$.
\endprclaim

\demo{Proof}
The proof proceeds by induction on $n+m$.
If $m = n = 2$, then it is easy to see that $\perm$ is prime, and of
the form (3).  So assume that $n+m\geq 5$.
Without loss of generality we may assume that $n\geq m$ and $n\geq 3$.  
Let $P$ be a prime ideal containing $\perm$.
We first show that $P$ contains all the entries from some row
or from some column of $M$.

Assume that there is no column all of whose entries
are in $P$. 
If $P$ contains all but one entry from each row of $M$, then,
since $m\leq n$,
there necessarily exists a $(2\times 2)$-submatrix 
$$
\pmatrix x_{ij}&x_{ik}\\x_{lj}&x_{lk}\endpmatrix
$$
of $M$ such that the entries of either the diagonal or codiagonal
do not lie in $P$, and the other two entries are elements of $P$,
say $x_{ij}$ and $x_{lk}$ do not lie in $P$.  
Since $x_{ij}x_{lk}+x_{lj}x_{ik}\in \perm\subset P$
and $x_{lj}x_{ik}\in P$, it follows that $x_{ij}x_{lk}\in P$, which is
a contradiction, since $P$ is prime.

So, necessarily, there is a row of $M$, which has two entries not
in $P$,  say $x_{11}, x_{12}\notin P$.  By Lemma~\twoone,
$x_{11}x_{12}x_{ij}\in \perm \subset P$ for all $i>1, j>2$.  
This implies that $x_{ij}\in P$ for all $i>1, j>2$.
By assumption every column of $M$ contains at least one entry not in $P$,
so that $x_{1j}\notin P$ for all $j>2$. 
Now consider $x_{21}$ and $x_{22}$.
Since $n\geq 3$, we know by Lemma~\twoone\ %
that $x_{21}x_{22}x_{13}\in P$.
Since $x_{13}\notin P$,
this implies that one of the other two factors, say $x_{21}$, is in $P$.
Since
$x_{11}x_{22}+x_{21}x_{12}\in P$, this implies that $x_{11}x_{22}\in P$.
Since $x_{11}\notin P$, we have $x_{22}\in P$.  Hence $P$ contains all
the entries from the second row of $M$.
Similarly,
one shows that it contains all the entries from rows $3,\ldots ,m$.
A similar argument applies if
$P$ does not contain all the entries from some row of $M$.

So we may assume that $P$ contains an ideal $J$
generated by all the entries from one row or one column of $M$.
By transposing $M$, if necessary,
we assume that $P$ contains a row of $M$.
Of course,
$P$ also contains the permanental ideal $P_2$ of the submatrix of $M$
obtained by deleting that row.
Let $Q$ be a prime ideal contained in $P$ which is minimal over $P_2$.
We have the containments
$$
\perm \subset P_2+J\subset Q+J\subset P,
$$
and $P$ is minimal over each of the smaller ideals.
We know 
the structure of $Q$ by induction hypothesis,
whence we know the structure of $Q + J$.
In particular,
we deduce that $Q + J$ is a prime ideal.
As $P$ is minimal over it,
$Q + J = P$.
If $Q$ is of types (1) or (3), respectively,
so is $P$.
Now suppose that $Q$ is of type (2).
Then $P$ is the ideal generated by all the entries in
$n-1$ columns of $M$, plus an extra entry.
But this extra entry makes $P$ not
minimal over $\perm$,
contradicting the assumption.
Thus $Q$ cannot be of type (2).
This completes the proof.  
\enddemo 

\remark{Remark}
Observe that the three types of primes in Theorem \minp\ %
satisfy no inclusion relations.
Primes of type (1) have height $(m-1)n$, those
of type (2) have height $m(n-1)$, and type (3) primes have height
$mn-4+1= mn-3$.
These heights are different in general.
Hence, $R/\perm$ is not equidimensional, that is, $\text{dim}(R/P)$ is
different for different minimal primes $P$.
This implies that
neither $\perm$ nor its radical is Cohen-Macaulay.   
In contrast, determinantal ideals are all radical \cite{BV, Corollary 5.8}
and Cohen-Macaulay \cite{BV, Theorem 5.3 and Corollary 5.17}.
\endremark

\prclaim{Corollary}
\label{\notCM}
If $(m,n) \not = (2,2)$ and
If $(m,n) \not = (3,3)$,
then $R/\perm$ is not equidimensional,
hence not Cohen-Macaulay.
\endprclaim

\prclaim{Corollary}
\label{\nocomps}
The number of minimal components of $\perm$ is
\roster
\item
$m+n+\binom m2\binom n2 \text{ if }m, n\geq 3$, 
\item
$m+\binom m2\binom n2 \text{ if }m\geq 3, n =2$, 
\item
$n+\binom m2\binom n2 \text{ if } m =2, n\geq 3$, 
\item
$1 \text{ if } m = n = 2$.
\endroster
\endprclaim

The explicit structure of the minimal primes over $\perm$
gives yet another property of
permanental ideals which contrasts with determinantal ideals.
Glassbrenner and Smith proved in \cite{GS}
that determinantal varieties have systems of
parameters that are highly symmetric and sparse.
We now show that this is not the 
case for $(2\times 2)$-permanental varieties.

By a parameter we mean
an element $a = \sum_{ij} a_{ij} x_{ij}$
with elements $a_{ij}$ in $R$
such that $a$ is not in any minimal prime of $\perm$.
A parameter is called {\it sparse}
if very few of the $a_{ij}$ are non-zero.
By Theorem~\minp,
for each $(2\times 2)$-submatrix of $M$,
one of the variables appearing in that submatrix
must have a non-zero coefficient in $a$,
for otherwise $a$ lies in a minimal prime of type (3).
Similarly,
if $n \ge 3$,
each row of $M$ contains an entry with a non-zero coefficient in $a$,
and if $m \ge 3$,
each column of $M$ contains an entry with a non-zero coefficient in $a$.
Thus no parameter on the $(2\times 2)$-permanental variety is sparse.

\section{Primary Decomposition.}

In this section we calculate a primary decomposition of $\perm$
and its radical.
Theorem~\minp\ determines all the minimal primes over $\perm$,
which are of three types.
We let ${\Cal P}_i, i =1, 2, 3$,
be the set of all minimal primes of type (i) as in Theorem~\minp.
Let ${\Cal P}={\Cal P}_1\cup{\Cal P}_2\cup{\Cal P}_3$.

\prclaim{Proposition}
\label{\minc}
The primary components of $\perm$ corresponding to the minimal
primes over $\perm$ are exactly the minimal primes in ${\Cal P}$
themselves.
\endprclaim

\demo{Proof}
The result is clear if $m = n = 2$.

Let $P \in {\Cal P}$.
Let $Q_P$ be the $P$-primary component of $\perm$.
Then $\perm\subset Q_P\subset P$, and $Q_P$ is characterized by
the property that $Q_P$ contains a power of $P$ and, 
if $rs\in Q_P$ and $r\notin P$, then $s\in Q_P$.
We will show that $Q_P$ contains all the generators of $P$, hence
is equal to $P$.

First let $P$ be a minimal prime over $\perm$
that is generated by all the entries of $M$ except for those in one row,
say the first one.
Then, by assumption, $n \ge 3$.
Let $x_{ij}, i>1$, be an entry of $M$.
Let $p\neq q\neq j \neq p$ be column labels
(which exist, since $n\geq 3$).
Then, by Lemma~\twoone, the element $x_{ij}x_{1p}x_{1q}$ is in 
$\perm\subset Q_P$.  Furthermore, $x_{1p}x_{1q}\notin P$, hence
$x_{ij}\in Q_P$.  This proves that $P = Q_P$.

The case where $P$ is generated by all entries of $M$ except for those
in one column is proved similarly.  

Finally, suppose that $P$ is
generated by one $(2\times 2)$-subpermanent and all entries outside
of this $(2\times 2)$-block, as in Part (3) of Theorem~\minp.
As above,
by using Lemma~\twoone\ %
and the defining property of $Q_P$ repeatedly,
all entries of $M$ outside of the $(2\times 2)$-block
are in $Q_P$.  Since any $(2\times 2)$-subpermanent is in 
$\perm\subset Q_P$, this completes the proof.
\enddemo

In order to compute the radical and a primary decomposition of $\perm$
we need to compute the intersection of all the minimal primes.
Define $I_i =\bigcap_{P\in{\Cal P}_i}P$.  
We first calculate $I_1$, $I_2$ and $I_3$,
which are the unmixed parts of $\perm$ of various dimensions.

\prclaim{Lemma}
\label{\Is}
Let $m, n\geq 2$.
Recall that $I_1$ is only defined if $n \ge 3$,
and $I_2$ is only defined if $m \ge 3$.
Then 
\roster
\item
$ I_1=\langle x_{ij}x_{lk}|i\neq l,  1\leq i, l\leq m,  
                1\leq j, k\leq n\rangle$;
\item
$ I_2=\langle x_{ij}x_{lk}|1\leq i, l\leq m,  1\leq j, k\leq n,  j\neq k
               \rangle$;
\item
$I_3=\perm +\langle \xipjqkr | \ijkm \rangle
            +\langle \xipjqkr | \pqrn \rangle$,
\endroster
where one or both of the first two ideals may be zero if $m$ or $n$ is $2$.
\endprclaim

\demo{Proof}
(1) Each $P\in {\Cal P}_1$ contains any product of entries from two
different rows,
which shows that the right-hand side is contained in $I_1$.
To show the other inclusion,
observe that each $P\in {\Cal P}_1$ is a monomial
ideal, and the intersection of monomial ideals is again monomial.
Let $u\in I_1$ be a monomial generator.  Then $u$ is a product of
entries of $M$.  Without loss of generality assume that $u = x_{11}v$
for some monomial $v$.  It is now sufficient to show that $v$ is
contained in the ideal generated by the entries of the last $m-1$
rows of $M$.  But this is clear, since $x_{11}v$ is contained in
all primes in ${\Cal P}_1$, in particular in the ideal that is 
generated by all entries except those in the first row.

The proof of (2) is similar to that of (1).
\enddemo

Our original proof of part (3) involved multiple lengthy 
and uninstructive induction arguments.
Instead we present Niermann's clever shortcut.
Then the proof of (3) follows by an easy application
of the following lemma from \cite{N}:

\prclaim{Lemma}
(Niermann \cite{N, p. 103})
\label{\niermann}
Let $R$ be an arbitrary ring
and $I_1, \ldots, I_l$,
$J_1, \ldots, J_l$ ideals in $R$
such that $I_i \subseteq J_j$ if $i \not = j$.
Then
$$
\bigcap_{i = 1} ^l (I_i + J_i)
= I_1 + \cdots + I_l + \bigcap_{i = 1} ^l J_i.
$$
\endprclaim

The proof, given in \cite{N}, is a straightforward induction on $l$.
This lemma finishes the proof of Lemma~\Is.

\medskip
We are now ready to compute the radical and primary decomposition
of $\perm$.

\prclaim{Theorem}
\label{\rad}
If $m$ or $n$ is equal to $2$,
then $\text{rad}(\perm)= \perm$.
If $m, n\geq 3$,  then
$$
\align
\text{rad}(\perm)&=I_1\cap I_2\cap I_3\\
&=\perm +\langle \xipjqkr |
i\neq j\neq k \neq i, p\neq q\neq r \neq p\rangle .
\endalign
$$
So $\perm$ is a radical ideal if and only if $m\leq 2$ or $n\leq 2$.
\endprclaim

\demo{Proof}
The result is of course trivial if $m = n = 2$.

Assume that $m = 2$ and $n \ge 3$.
Then $\perm$ has minimal primes of types (1) and (3) only
(cf.\ Theorem~\minp).
Thus
$$
\align
\text{rad}(\perm)&=I_1\cap I_3\\
&= \langle x_{ij}x_{lk}|i\neq l,  1\leq i, l\leq m,  
                1\leq j, k\leq n\rangle \\
&\qquad \bigcap
\left(\perm +\langle \xipjqkr | 1 \le i < j < k \le 2
\text{ or } \pqrn \rangle\right)\\
&= \perm +
\langle x_{ij}x_{lk}|i\neq l  \rangle 
\cap \langle \xipjqkr | \pqrn \rangle, \\
\endalign
$$
and the last intersection is in $\perm$ by Lemma~\twoone\ %
as $m = 2$.

Similarly,
if $n = 2$,
$\perm$ is a radical ideal.

Now assume that $m, n \ge 3$.
Then
$$
\align
\text{rad}(\perm)&=I_1\cap I_2 \cap I_3\\
&= 
\langle x_{ij}x_{lk}|j\neq k, i\neq l\rangle \cap I_3 \\
&= 
\langle x_{ij}x_{lk}|j\neq k, i\neq l\rangle \\
&\qquad \cap
\left(\perm +\langle \xipjqkr |
\dijkm \text{ or } \dpqrn \rangle\right)\\
&= \perm +
\langle \xipjqkr | \dijkm \text{ and } \dpqrn \rangle, \\
\endalign
$$
the latter equality by Lemma~\twoone.
Now, the monomials of the form $\xipjqkr$ with
distinct $i,j,k$ and distinct $p,q,r$ are not in $\perm$,
as can be easily verified via our Gr\"obner basis in Theorem~\gb.
This finishes the proof of the theorem.
\enddemo

Thus, for $m, n \ge 3$,
$P_2(M)$ is an example of an ideal
whose radical requires generators of degree higher than those of the
ideal itself.
Thus permanental ideals might seem a possible candidate
for a negative answer to a question of Ravi \cite{R},
whether for a homogeneous ideal
its Castelnuovo-Mumford regularity is at least as big as the
regularity of its radical.
This is not the case, however.

\prclaim{Corollary}
\label{\reg}
Let $I$ be the radical of $\perm$.  Then
\roster
\item
$\text{reg} (\perm)\geq \text{reg}(I)$;
\item
$\text{reg}(\perm)\geq 1+\sum_{i\geq3}{m\choose i}{n\choose i}$.
\endroster
\endprclaim

\demo{Proof}
Consider the short exact sequence of graded $R$-modules 
$$
0\longrightarrow I/\perm \longrightarrow R/\perm \longrightarrow R/I\longrightarrow 0.
$$
Since $I/\perm$ has finite length by Theorem~\thmfinite,
it follows from \cite{E, Corollary 20.19.d}
that 
$$
\text{reg}(R/\perm)=\text{max}\{\text{reg}(I/\perm, R/I\}\geq \text{reg}(R/I).
$$
As for any ideal $J \subset R$,
$\text{reg}(J)=\text{reg}(R/J)+1$,
the first assertion of the corollary follows.

The displayed formula also shows that
$\text{reg}(R/\perm) \geq \text{reg}(I/\perm)$,
and as $I/\perm$ has finite length,
$\text{reg}(I/\perm) = \text{length}(I/\perm)$.
Thus
$$
\text{reg}(\perm)= 1+\text{reg}(R/\perm)\ge 1+\text{length}(I/\perm)=
1+\sum_{i\geq 3}{m\choose i}{n\choose i}.
$$
Here, the first inequality follows from the fact that $I/\perm$ has finite 
length, and the second equality follows from Theorem 3.3.
\enddemo

\prclaim{Corollary}
\label{\emb}
$\perm$ has embedded components if and only if $m, n \ge 3$.
\endprclaim

\ \

It follows that for 
all $M$ with $m, n \ge 3$,
$\perm$ has embedded components.
We now show that,
in fact, $\perm$ has exactly one embedded component.

Let $Q$ be the ideal
$$
Q=\perm +\langle x_{ij}^2|1\leq i\leq m,1\leq j\leq n\rangle .
$$
Then for $m, n\geq 3$,  we obtain from Theorem~\rad\ %
that
$$
\align
Q\cap I_1\cap I_2\cap I_3
&=\left(\perm +\langle x_{ij}^2\rangle\right)\cap
\left(\perm +\sum\langle x_{ij}x_{kl}x_{pq}\rangle\right)\\
&=\perm +
\langle x_{ij}^2\rangle\cap
\left(\perm +\sum\langle x_{ij}x_{kl}x_{pq}\rangle\right)\\
\endalign
$$
The sum in this identity extends
over all subscripts $i\neq k\neq p\neq i$
and $j\neq l\neq q\neq j$.
To simplify the last intersection and prove that it lies in $\perm$,
we use Gr\"obner bases again, in any diagonal term order.
It suffices to show that for any element $f$ in the intersection,
after reducing $f$ with respect to the Gr\"obner basis of $\perm$
as in Theorem~\gb,
we get $f = 0$.
If this is false,
the leading monomial $x$ of $f$ is non-zero,
and hence necessarily a multiple of one of the monomials of type (6)
in the Gr\"obner basis for
$I_1\cap I_2\cap I_3$
$= \sqrt{\perm}$,
as was computed in Theorem~\gbrad.
Moreover,
$x$ has to be divisible by the square of a variable $x_{ij}$.
But least common multiples of these two types of
monomials are all in the ideal generated by
monomials of types (2)-(6) in the Gr\"obner basis of $\perm$
(as in Theorem~\gb),
contradicting the initial statement that $f$ is reduced
with respect to this basis.

Since $Q$ is primary to the homogeneous maximal ideal of all variables
and each of $I_1,I_2,I_3$ is the intersection of distinct minimal
components of $\perm$,
we have calculated a primary decomposition of $\perm$.

\prclaim{Theorem}
\label{\pd}
Let $m, n\geq 3$.  The intersection
$$
\perm =Q\cap I_1\cap I_2\cap I_3,
$$
after rewriting each $I_i$ as the intersection
of the minimal primes of type (i),
is an irredundant primary decomposition of $\perm$.
\endprclaim

\vfill\eject

\enddocument